 \documentclass[10pt]{amsart}

   \newtheorem{lemma}{Lemma}[section]
   \newtheorem{theorem}[lemma]{Theorem}

   \newtheorem{definition}[lemma]{Definition}

\title[]{}

\author{Jinqiao Duan}
\address[J. Duan]
{Department of Applied Mathematics\\
 Illinois Institute of Technology\\
   Chicago, IL 60616, USA \\and\\Department of Mathematics\\
    The University of Science and Technology of China\\ Hefei 230026, China}
\email[J.~Duan]{duan@iit.edu}
\author{Kening Lu}
\address[K. Lu]
{Department of Mathematics\\
    Brigham Young University\\
    Provo, Utah 84602, USA \\and \\ Department of Mathematics\\ Michigan State University \\East Lansing,
    Michigan 48824, USA }
\email[K.~Lu]{klu@math.byu.edu or klu@math.msu.edu}

\author{Bj{\"o}rn Schmalfu{\ss}}
\address[B. Schmalfu{\ss}]{%
  Department of Sciences\\
  University of Applied Sciences\\
  Geusaer Stra{\ss}e\\
  06217 Merseburg, Germany}
\email[B.~Schmalfu{\ss}]{bjoern.schmalfuss@in.fh-merseburg.de}

\title[Stable and unstable  manifolds for  stochastic PDEs]
{Smooth stable and unstable  manifolds for stochastic partial
 differential equations}

\dedicatory{Dedicated to Professor Shui-Nee Chow on his  $60^{th}$
birthday}

\keywords{Invariant manifolds, cocycles, non-autonomous dynamical
systems, stochastic partial differential equations, generalized
fixed points. \\ This work was partially supported by NSF 0209326,
NSF0200961, and a travel grant from German DFG Schwerpunktprogramm
{\em Interagierende zuf\"allige Systeme von hoher Komplexit\"at}.}

\begin{document}

\begin{abstract}
{\bf J. Dynamics and Diff. Eqns 2004, in press.}

 Invariant manifolds are fundamental tools for describing and
understanding nonlinear  dynamics. In this paper, we present a
theory of stable and unstable  manifolds for infinite dimensional
random dynamical systems generated by
 a class of   stochastic  partial differential equations.
We  first show the existence of Lipschitz continuous stable and
unstable manifolds by the Lyapunov-Perron's method. Then, we prove
the smoothness of these invariant manifolds.

\end{abstract}

\maketitle

%%%%   Section 1   %%%%%%%%%%
%%%%%%%%%%%%%%%%%%%%%
\section{\bf Introduction}
This paper, which is a sequel to \cite{DLS1}, is devoted to the
existence and smoothness of stable and unstable manifolds for a
class of stochastic partial differential equations (PDEs).\\

We consider a nonlinear stochastic evolution equation with a
multiplicative white noise:
\begin{equation} \label{eq(1.1)}
\frac{du}{dt}=A u+F(u)+u\,\dot{W},
\end{equation}
where $A$ is a generator of a $C_0$-semigroup $e^{At}$ satisfying
an exponential dichotomy condition, $F(\phi)$ is a Lipschitz
continuous operator with $F(0)=0$, and $u\,\dot{W}$ is a noise.
The precise conditions on them will given in the next section.
Some physical systems or fluid systems with noisy perturbations
proportional to the state of the system may be modeled
by this equation.  \\

In \cite{DLS1}, we proved the existence of Lipschitz continuous
unstable manifolds for stochastic partial differential equation
(\ref{eq(1.1)}) by using a random graph transform and a
generalized fixed point theorem.\\

In the present paper, we study the existence and smoothness of
stable and unstable manifolds for equation (\ref{eq(1.1)}). In
brief, our main results on the stable and unstable manifolds
manifolds may be summarized as follows (the precise statements are
given in Sections 3-5). We assume that the semigroup $e^{At}$
satisfies an exponential condition and the Lipschiz constant of
$F$ is dominated by a spectrum gap. Then, there exist global
Lipschitz continuous stable and unstable manifolds for stochastic
partial differential equation (\ref{eq(1.1)}). Furthermore, if $F$
is $C^k$ and a large spectrum gap condition holds, then these
stable and unstable manifolds are $C^k$ smooth. The manifolds we
study here actually are so-called pseudo-stable and
pseudo-unstable manifolds which include the usual stable and
unstable manifolds. As for the deterministic systems, we do not
need the large spectrum condition for the smoothness of the usual
stable and unstable manifolds of stochastic partial differential
equation (\ref{eq(1.1)}). When $F$ is a $C^1$ function, the
conditions for existence of the $C^1$ stable and unstable
manifolds are the same
as those for the existence of Lipschitz manifolds.\\

In this paper, we also consider a nonlinear stochastic evolution
equation with a additive white noise:
\begin{equation} \label{eq(1.1)'}
\frac{du}{dt}=A u+F(u)+ \dot{W}.
\end{equation}
The precise conditions on them will given in the next section. We
will see that after random transformations, equation
(\ref{eq(1.1)'}) can be regarded as a
special case of equation (\ref{eq(1.1)}).\\

The approach we use here is based on the Lyapunov-Perron's method.
This approach differs from the Hadamard's graph transform
method that we used in \cite{DLS1}.
However, the optimal conditions for the existence of Lipschitz
unstable manifolds obtained by these two different methods are the same. \\

Stable and unstable manifolds play  an important role in the study of
 nonlinear dynamical systems. Hadamard \cite{Had01} constructed
the unstable manifold of a hyperbolic fixed point of a
diffeomorphism of the plane using a geometric method. This
geometric method is now called Hadamard's graph transform.
Lyapunov \cite{Lia47} and Perron \cite{Per28} constructed the
unstable manifold for an equilbrium point by formulating the
problem in terms of an integral equation. This method is analytic
rather than geometric and now is called the method of  Lyapunov
and Perron. There is an extensive literature on stable and
unstable manifolds for both finite and infinite dimensional
deterministic dynamical systems; see Henry \cite{Hen81},  Babin
and Vishik \cite{BabVis89} or Bates et. al. \cite{BatLuZen98} and
the references therein.\\

Recently, there are some works on invariant manifolds for
stochastic  ordinary differential equations by Wanner
\cite{Wanner}, Arnold \cite{Arn98}, Mohammed and Scheutzow
\cite{MohScheu99}, and Schmalfu{\ss} \cite{Schm97a}. Wanner's
method is based  on the Banach fixed point theorem on a
space composed of functions with particular exponential growth
rates. This method is   essentially   the Liapunov-Perron approach.
In contrast to this
method, Mohammed and Scheutzow have applied a classical technique
due to Ruelle \cite{Rue82}  to the stochastic  ordinary differential
equations driven by semimartingals.  Caraballo et. al.
\cite{CarLanRob01} have considered   invariant manifolds for a
stochastic reaction diffusion  equation. \\

In Section 2, we recall some  basic concepts and results for
random dynamical systems and   stochastic partial differential
equations (PDEs). We then prove the existence of the Lipschitz
stable manifold for the   stochastic  PDE  (\ref{eq(1.1)}) in
Section 3. In Section 4, we prove the smoothness of the stable
manifold. The results on the unstable manifold for (\ref{eq(1.1)})
are given in Section 5.

%%%% Section 2 %%.
%%%%%%%%%%%%%%%%%%%%
\section{\bf Stochastic PDEs and Random Dynamical Systems}

In this section, we introduce some basic notations, assumptions,
concepts, and results on stochastic partial differential equations
and random dynamical systems.

\subsection{Stochastic PDEs with a Multiplicative Noise}
Let $H$ be an infinite dimensional  separable  Hilbert  space with
norm $|\cdot|$. Consider the nonlinear stochastic partial
differential equation
\begin{equation}\label{eq(2.1)}
\frac{d u}{dt}=Au+F(u)+u\,\dot{W},
\end{equation}
where $u\in H$, $W(t)$ is the standard $\mathbb{R}-$valued Wiener
process on a probability space $(\Omega,\mathcal{F},\mathbb{P})$,
and the generalized time-derivative $\dot{W}$ formally describes a
{\em white-in-time  noise}. Note that $u \,\dot{W}$ is interpreted
as a Stratonovich stochastic differential.
\\

We assume that the linear operator $A: D(A)\to H$ generates a
strongly continuous semigroup $e^{At}$ on $H$, which satisfies the
exponential dichotomy with exponents $\alpha > \beta$ and bound
$K$, i.e., there exists a continuous projection $P^+$ on $H$ such
that
\begin{itemize}
\item[(i)] $P^+e^{At}=e^{At}P^+$;
\item[(ii)] the restriction $e^{At}|_{R(P^+)}$, $t\geq 0$, is an
isomorphism of the range ${R(P^+)}$ of $P^+$ onto itself, and we
define $e^{At}$ for $t <0$ as the inverse map.
\item[(iii)]
\begin{align}\label{dichotomy}
\begin{split}
&|e^{At}P^+ x|\le K e^{\alpha t} |x|,\quad t\le 0, \\
&|e^{At} P^-x|\le K e^{\beta t} |x|,\quad t\ge 0,\;
\end{split}
\end{align}
\end{itemize}
where $P^-=I-P^+$. Denote $H^-= P^-H$ and $H^+=P^+H$. Then,
$H=H^+\oplus H^-$. We will call $H^-$ and $H^+$ the stable
subspace  and the unstable subspace, respectively.

\vskip0.1in The nonlinear term $F$ satisfies $F(0)=0$ and is assumed to be Lipschitz
continuous on $H$
\[
\|P^\pm (F(x_1)-F(x_2))\|_H\le \hbox{{\rm Lip}} F \|x_1-x_2\|_H
\]
with the Lipschitz constant $\hbox{{\rm Lip}} F>0$.

\vskip0.1in

The existence theory for stochastic evolution equations is usually
formulated for Ito equations as in Da Prato and Zabczyk
\cite{DaPZab92}, Chapter 7. The equivalent Ito equation for
(\ref{eq(2.1)}) is given by
\begin{equation}\label{eq(2.2)}
du=Au\,dt+F(u)\,dt+\frac{u}{2}dt+u\,dW.
\end{equation}
Then, for any initial data $x\in H$, there exists a unique
solution of (\ref{eq(2.2)}). For details about the properties of
this solution see Da Prato and Zabczyk \cite{DaPZab92}, Chapter 7.
\\

The solution of (\ref{eq(2.2)}) can be written as a {\em mild}
solution:
\[
u(t)=e^{At}x+\int_0^t(e^{A(t-s)}(F(u(s))+\frac{u(s)}{2})ds
+\int_0^te^{A(t-s)}u(s)dW,\quad x\in H
\]
almost surely for any $x\in H$. Note that the theory in
\cite{DaPZab92}
requires that the associated probability space
 $(\Omega,\mathcal{F},\mathbb{P})$ is  complete.

\subsection{Random Dynamical Systems}

Let us first look at flows on the probability space
$(\Omega,\mathcal{F},\mathbb{P})$. A
{\em flow} $\theta$ of mappings $\{\theta_t\}_{t\in\mathbb{R}}$ is
defined  on the sample space $\Omega$   such that
\begin{equation}\label{eq(2.3)}
\theta:\mathbb{R}\times \Omega\to \Omega,\qquad \theta_0={\rm
id}_\Omega, \qquad\theta_{t_1}\circ\theta_{t_2}= \theta_{t_1+t_2}
\end{equation}
for $t_1,\,t_2\in\mathbb{R}$. This flow is supposed to be
$(\mathcal{B}(\mathbb{R})\otimes\mathcal{F},\mathcal{F})$-measurable,
where $\mathcal{B}(\mathbb{R})$ is the $\sigma$-algebra of Borel
sets on the real line $ \mathbb{R}$. To have this measurability,
it is not allowed to replace $\mathcal{F}$ by its
$\mathbb{P}$-completion $\mathcal{F}^\mathbb{P}$; see Arnold
\cite{Arn98} Page 547. In addition, the measure $\mathbb{P}$ is
assumed to be ergodic with respect to
$\{\theta_t\}_{t\in\mathbb{R}}$. Then
$\theta:=(\Omega,\mathcal{F},\mathbb{P},\mathbb{R},\theta)$ is
called a
 metric dynamical system.\\

For the SPDE's with a multiplicative noise, we will consider a
special but very important metric dynamical system  induced by the
Wiener process. Let $W(t)$ be a two-sided Wiener process with
trajectories in the space $C_0(\mathbb{R},\mathbb{R})$ of real
continuous functions defined on $\mathbb{R}$,
  taking zero value  at $t=0$. This set is
equipped with the compact open topology. On this set we consider
the measurable flow $\theta=\{\theta_t\}_{t\in\mathbb{R}},\,
\mbox{defined by} \,\theta_t\omega=\omega(\cdot+t)-\omega(t)$. The
distribution of this process is a measure on
$\mathcal{B}(C_0(\mathbb{R},\mathbb{R}))$ which is called the {\em
Wiener measure}. Note that this measure is ergodic with respect to
the above flow; see the Appendix in Arnold \cite{Arn98}. Later on
we will consider, instead of the whole
$C_0(\mathbb{R},\mathbb{R})$,
 a
$\{\theta_t\}_{t\in\mathbb{R}}$-invariant subset $\Omega\subset
C_0(\mathbb{R},\mathbb{R})$ of $\mathbb{P}$-measure one and the
trace $\sigma$-algebra $\mathcal{F}$ of
$\mathcal{B}(C_0(\mathbb{R},\mathbb{R}))$ with respect to
$\Omega$. A set $\Omega$ is called
$\{\theta_t\}_{t\in\mathbb{R}}$-invariant if
$\theta_t\Omega=\Omega$ for $t\in\mathbb{R}$. On $\mathcal{F}$ we
consider the restriction of the
Wiener measure also denoted by $\mathbb{P}$.\\

The dynamics of the system on the state space $H$ over the {\em driven}
 flow
$\theta$ is described by a cocycle. For our applications it is
sufficient to assume that $(H,d_H)$ is a complete metric space. A
cocycle $\phi$ is a mapping:
\[
\phi:\mathbb{R}^+\times \Omega\times H\to H
\]
which  is
$(\mathcal{B}(\mathbb{R})\otimes\mathcal{F}\otimes\mathcal{B}(H),\mathcal{F})$-measurable
 such that
\begin{equation*}
%\label{eq3}
\begin{split}
&\phi(0,\omega,x)=x \in H,\\
&
\phi(t_1+t_2,\omega,x)=\phi(t_2,\theta_{t_1}\omega,\phi(t_1,\omega,x)),
\end{split}
\end{equation*}
for $t_1,\,t_2\in\mathbb{R}^+,\,\omega\in \Omega,$ and $x\in H$.
Then $\phi$ together with the metric dynamical system $\theta$
forms a {\em random dynamical system}.

\subsection{Conjugated Random PDEs}

\vskip0.1in In \cite{DLS1}, we used a coordinate transform to
convert conjugately a stochastic partial differential equation
into an infinite dimensional random dynamical system. Although it
is well-known that a large class of partial differential equations
with stationary random coefficients as well as  Ito stochastic ordinary
differential equations generate random dynamical systems (for
details see Arnold \cite{Arn98}, Chapter 1), this problem is still
unsolved for stochastic partial differential equations with a
general noise term $C(u)\,dW$. The reasons are: (i) The stochastic
integral is only defined almost surely where the exceptional set
may depend on the initial state $x$;     (ii)
  Kolmogorov's theorem,   as cited in
Kunita \cite{Kun90} Theorem 1.4.1,  is only true for finite
dimensional random fields; and (iii)    the cocycle has to be
defined for  {\em any} $\omega\in\Omega$.  Nevertheless, for the
noise term $u\,dW$ considered here,  we can show  that  the
stochastic PDE (\ref{eq(2.2)}) indeed generates a random dynamical
system.

 \vskip0.1in\noindent We considered a   linear stochastic
differential equation:
\begin{equation}\label{eq(2.4)}
dz+z\,dt=dW.
\end{equation}
A solution of this equation is called an Ornstein-Uhlenbeck
process. We have the following result.
\begin{lemma}\label{lem(2.1)}
i) There exists a $\{\theta_t\}_{t\in\mathbb{R}}$-invariant set
$\Omega\in\mathcal{B}(C_0(\mathbb{R},\mathbb{R}))$ of full measure
with sublinear growth:
\[
\lim_{t\to\pm\infty}\frac{|\omega(t)|}{|t|}=0,\qquad
\omega\in\Omega
\]
of $\mathbb{P}$-measure one.\\
ii) For $\omega\in\Omega$ the random variable
\[
z(\omega)=-\int_{-\infty}^0e^\tau\omega(\tau)d\tau
\]
exists and generates a unique stationary solution of (\ref{eq(2.4)})
given by
\[
\Omega\times \mathbb{R}\ni(\omega,t)\to z(\theta_t\omega)
=-\int_{-\infty}^0e^\tau\theta_t\omega(\tau)d\tau
=-\int_{-\infty}^0e^\tau\omega(\tau+t)d\tau+\omega(t).
\]
The mapping $t\to z(\theta_t\omega)$
is continuous.\\
iii) In particular,
\begin{equation*}
%\label{eq401}
\lim_{t\to\pm\infty}\frac{|z(\theta_t\omega)|}{|t|}=0\quad
   \text{for }\omega\in\Omega.
\end{equation*}
iv) In addition,
\begin{equation*}
\lim_{t\to\pm\infty}\frac{1}{t}\int_0^tz(\theta_\tau\omega)d\tau=0\quad
  \text{for }\omega\in\Omega.
 \end{equation*}
\end{lemma}

We now replace $\mathcal{B}(C_0(\mathbb{R},\mathbb{R}))$ by
\[
\mathcal{F}=\{\Omega\cap F,\,F\in
\mathcal{B}(C_0(\mathbb{R},\mathbb{R}))\}
\]
for $\Omega$ given in Lemma \ref{lem(2.1)}. The probability
measure is the restriction of the Wiener measure to this new
$\sigma$-algebra, which is  also denoted by $\mathbb{P}$. In the
following we will consider the metric dynamical system
\[
(\Omega,\mathcal{F},\mathbb{P},\mathbb{R},\theta).
\]
We show that the solution of (\ref{eq(2.2)}) defines a random
dynamical systems. To see this we
 consider the following partial differential equation with random coefficients
\begin{equation}\label{eq(2.5)}
\frac{du}{dt}=Au+z(\theta_t\omega)u+G(\theta_t\omega,u),\quad
u(0)=x\in H
\end{equation}
where $G(\omega,u):=e^{z(\omega)}F(e^{-z(\omega)}u)$. It is easy
ro see that for any $\omega\in\Omega$ the function $G$ has the
same global Lipschitz constant $L$ as $F$. In contrast to the
stochastic PDE  (\ref{eq(2.2)}),  no stochastic differential appears in the random PDE  (\ref{eq(2.5)}).
The solution can be interpreted in a mild sense
\begin{equation}\label{eq(2.6)}
u(t)=e^{At +\int_0^tz(\theta_\tau\omega)d\tau}x+\int_0^t
e^{A(t-s)+\int_s^tz(\theta_r\omega)dr} G(\theta_s\omega,u(s))ds.
\end{equation}
We note that this equation has  a unique solution for each
$\omega\in\Omega$. No exceptional sets appear.
 Hence the solution mapping
\[
(t,\omega,x)\to u(t,\omega, x)
\]
generates a random dynamical system. Indeed, the mapping $u$ is
$(\mathcal{B}(\mathbb{R})\otimes\mathcal{F}\otimes\mathcal{B}(H),\mathcal{F})$-measurable.\\

%Let now $\hat u(t,\omega,x)$ be the solution mapping of
%(\ref{eq(2.2)}) which is defined for
%$\omega\in\Omega\in\mathcal{F}^\mathbb{P},\,\mathbb{P}(\Omega)=1$.
We now introduce the transform
\begin{equation}\label{eq(2.7)}
T(\omega,x)=xe^{-z(\omega)}
\end{equation}
and its inverse transform
\begin{equation}\label{eq(2.8)}
T^{-1}(\omega,x)=xe^{z(\omega)}
\end{equation}
for $x\in H$ and $\omega\in\Omega$.
\begin{lemma}\label{lem(2.2)}
Suppose that $u$ is the random dynamical system generated by
(\ref{eq(2.5)}) Then
\begin{equation}\label{eq(2.9)}
(t,\omega,x)\to T^{-1}(\theta_t\omega,
u(t,\omega,T(\omega,x)))=:\hat u(t,\omega,x)
\end{equation}
is a random dynamical system. For any $x\in H$ this process
$(t,\omega)\to\hat u(t,\omega,x)$ is a solution to
(\ref{eq(1.1)}).
\end{lemma}

Similar transformations were used in Caraballo, Langa and
Robinson \cite{CarLanRob01}, and Schmalfu{\ss} \cite{Schm97c}.
%Note
%that our transform has the  advantage that the solution of
%(\ref{eq(2.5)}) generates a random dynamical system for the
%$\omega$-wise differential equation in \cite{CarLanRob01} which is
%based on the transform
%$e^{W(t)}$.\\

\subsection{Stochastic PDEs with a Additive white noise.}
 We mention another application. We consider a stochastic
evolution equation with an additive white noise
\begin{equation}\label{add1}
    \frac{d\hat u}{dt}=A\hat u+F(\hat u)+\dot{W},\qquad \hat
    u(0)=x
\end{equation}
where $\dot{W}$ is a white noise given as the generalized temporal
derivative of a Wiener process with continuous paths in $H$. For
simplicity we suppose that $W$ has a covariance with finite trace.
For a comprehensive presentation of these equations see
\cite{DaPZab92}. For this problem we have to choose a similar
metric dynamical system as above but $\Omega$ is contained in the
space of trajectories $C_0(\mathbb{R},H)$.
\\
Suppose that $u^\ast$ is a stationary solution to (\ref{add1}).
This means that for the random variable $u^\ast$  with values in
$H$ defined on a  $\{\theta_t\}_{t\in\mathbb{R}}$-invariant set of
full measure
\begin{equation*}
    t\to u^\ast(\theta_t\omega)
\end{equation*}
is a solution version for (\ref{add1}). It will not be the topic
of this article to deal with  stationary solutions. For the
existence of stable stationary solutions see Caraballo et al.
\cite{CarKloSchm03}.\\
We now define the nonlinear operator
\begin{equation*}
    G(\omega,x)=F(x+u^\ast(\omega))-F(u^\ast(\omega)).
\end{equation*}
Note that $G$ has the same Lipschitz constant as $F$. In addition,
$G(\omega,0)=0$. Hence, the problem
\begin{equation}\label{add2}
    \frac{du}{dt}=Au+G(\theta_t\omega,u),\qquad u(0)=x\in H.
\end{equation}
has a stationary solution which is identical zero. We introduce
the random transformations
\begin{equation*}
    T(\omega,x)=x-u^\ast(\omega),\qquad
    T^{-1}(\omega,x)=x+u^\ast(\omega).
\end{equation*}
\begin{lemma}\label{lem(2.2A)}
Suppose that $u$ is the random dynamical system generated by
(\ref{add2}). Then
\begin{equation*}
    T^{-1}(\theta_t\omega,u(t,\omega,
    T(\omega,x)))=:\hat u(t,\omega,x)
\end{equation*}
is a random dynamical system. For any $x\in H$ the process
\begin{equation*}
    (t,\omega)\to \hat u(t,\omega,x)
\end{equation*}
is a solution for (\ref{add1}).
\end{lemma}

We notice that equation (\ref{add2}) can be regarded as equation
(\ref{eq(2.5)}) with $z=0$. We refer to \cite{CaDuLuSchm03} for
more general cases. For the remainder of this article, we consider
only equation (\ref{eq(2.5)}).\\

\subsection{Definition of Invariant Manifolds}

We first recall that a multifunction
$M=\{M(\omega)\}_{\omega\in\Omega}$ of nonempty closed sets
$M(\omega),\,\omega\in\Omega$, contained in a complete separable
metric space $(H,d_H)$ is called a {\em random set} if
\[
\omega\to\inf_{y\in M(\omega)}d_H(x,y)
\]
is a random variable for any $x\in H$.
\begin{definition}
A random set $M(\omega)$ is called an invariant set for a random
dynamical system $\phi(t, \omega, x)$ if we have
\[
\phi(t,\omega,M(\omega))\subset M(\theta_t\omega)\; \hbox{for} \;
t\geq 0.
\]
If we can represent $M$ by a graph of a $C^k$ (or Lipschitz)
mapping
\[
h^s(\cdot,\omega): H^-\to H^+
\]
such that
\begin{align*}
M(\omega)=M^s(\omega)=\{\xi + h^s(\xi,\omega) | \xi\in H^-\}
\end{align*}
then $M^s(\omega)$ is called a $C^k$ (or Lipschitz) stable
manifold, where $H^-$ is the stable subspace and $H^+$ is the
unstable subspace, which are introduced in Section 2.1.

 If we can represent $M$
by a graph of a $C^k$ (or Lipschitz) mapping
\[
h^u(\cdot, \omega): H^+\to H^-
\]
such that
\begin{align*}
M(\omega)=M^u(\omega)=\{\xi + h^u(\xi, \omega)|\xi \in H^+\}
\end{align*}
then $M^u(\omega)$ is called a $C^k$ (or Lipschitz) unstable
manifold.
\end{definition}

%%%%%%%%%%%%%
%%%%%%%%%%%%%
%%%% Section 3 %%%%
%%%%%%%%%%%%%%
\section{\bf Lipschitz Stable Manifolds}

In this section, we first show the existence of a Lipschitz
continuous stable manifold for the random partial differential
equation
\begin{equation}\label{eq(3.1)}
\frac{d u}{dt}=Au+z(\theta_t\omega)u+G(\theta_t\omega,u),\quad
u(0)=u_0\in H
\end{equation}
Then, we apply the inverse transformation $T^{-1}$ to get a stable
manifold for the stochastic partial differential equation
(\ref{eq(2.2)}).\\

Denote by $u(t, \omega, u_0)$ the solution of  (\ref{eq(3.1)}) with the
initial data $u (0, \omega, u_0) = u_0$. We define the Banach
Space for each $\eta$, $\beta<\eta<\alpha$
\[
C_\eta^+=\{\phi:[0, \infty)\to H \;|\;\phi \hbox{  is continuous
and }\sup_{t\in [0, \infty)}e^{-\eta t - \int_0^t
z(\theta_\tau\omega) d\tau} |\phi(t)| < \infty\}
\]
with the norm
\[|\phi|_{C_\eta^+}=\sup_{t\in [0, \infty)}e^{-\eta t-\int_0^t z(\theta_\tau\omega)
d\tau} |\phi(t)|.
\]

Let
\[M^s(\omega)=\{u_0\in H\;|\; u(\cdot, u_0, \omega) \in C_\eta^+\}
\]
This is the set of all initial datum through  which solutions
decay  as $e^{\eta t +\int_0^t z(\theta_\tau\omega) d\tau}$.
  We
shall prove that $M^s(\omega)$ is invariant and is given by the
graph of a Lipschitz function.

\begin{theorem} \label{Thm(3.1)}    If $$K \;  \hbox{{\rm Lip}}_u G \; (\frac{1}{\eta-\beta}+\frac{1}{\alpha-\eta})<
1,$$ then there exists a Lipschitz invariant  stable manifold for
 the random partial differential  equation (\ref{eq(3.1)}) which is given by
\[
M^s(\omega) = \{\xi+h^s(\xi)\big | \xi \in H^-\},
\]
where $h^s : H^-\to  H^+$ is a Lipschitz continuous mapping and
satisfies $h^s(0)=0$. Note that $k, \alpha, \beta$ are from the
exponential dichotomy condition (\ref{dichotomy}) and $\hbox{{\rm
Lip}}_u G$ denotes the Lipschitz constant of $G(\cdot, u)$ with
respect to $u$.
\end{theorem}

\vskip0.1in \noindent {\bf Remark}: $\eta=(\alpha+\beta)/2$
minimizes the quantity
$$
K \; \hbox{{\rm Lip}}_u G \;
(\frac{1}{\eta-\beta}+\frac{1}{\alpha-\eta})
$$

\begin{proof}
We will  show that $M^s(\omega)$ is given by the graph of a
Lipschitz function over $H^-$.  First we claim that $u^0 \in \
M^s(\omega)$ if and only if  there exists a function $u(\cdot)\in
C_{\eta}^+$ with $u(0) = u^0$ and satisfies
\begin{align}\label{eq(3.2)}\begin{split}
u(t)=&e^{At+\int_0^tz(\theta_s\omega)ds}\xi+\int_0^t
e^{A(t-s)+\int_s^tz(\theta_r\omega)dr}
P^-G(\theta_s\omega,u(s))ds\\
&+\int^t_\infty e^{A(t-s)+\int_s^tz(\theta_r\omega)dr}
P^+G(\theta_s\omega,u(s))ds.
\end{split}
\end{align}
where $\xi = P^- u^0$.

To prove this claim,  first we let $u^0 \in M^s(\omega)$.  By
using the variation of constants formula, we have  that
\begin{align}\label{eq(3.3)}\begin{split}
P^-u(t,u^0, \omega) = &e^{A t + \int_0^tz(\theta_s\omega)ds}
P^-u^0 \\ &+ \int_0^te^{A(t-s)+\int_s^tz(\theta_r\omega)dr}
P^-G(\theta_s\omega,u)ds.
\end{split}
\end{align}

\begin{align}\label{eq(3.4)}\begin{split}
P^+ u(t,u^0, \omega) =&e^{A(t-\tau)+
\int_\tau^tz(\theta_s\omega)ds}P^+u(\tau,u^0, \omega)\\
&+\int_{\tau}^t e^{A(t-s)+\int_s^tz(\theta_r\omega)dr}
P^-G(\theta_s\omega,u)ds.
\end{split}
\end{align}
Since $u \in C_\eta^+$, we have  for $t < \tau , 0 < \tau$ that

\begin{align*}
&|e^{A(t-\tau)+ \int_\tau^tz(\theta_s\omega)ds}P^+u(\tau,u^0,
\omega)| \\
&\leq e^{\alpha(t-\tau)} e^{\int_0^tz(\theta_s\omega)ds} e^{\eta
\tau} |u|_{C_\eta^+} \\
& = e^{\alpha t +\int_0^tz(\theta_s\omega)ds}
e^{-(\alpha-\eta)\tau} \to 0\;\hbox{as }\tau \to +\infty.\\
\end{align*}

Then, taking the limit $\tau \to + \infty$ in (\ref{eq(3.4)}), we
have that
\begin{align}\label{eq(3.5)}
P^+u(t,u^0, \omega)=\int^t_\infty
e^{A(t-s)+\int_s^tz(\theta_r\omega)dr}
P^+G(\theta_s\omega,u(s))ds.
\end{align}

Combining  (\ref{eq(3.3)}) and (\ref{eq(3.5)}), we have
(\ref{eq(3.2)}). The converse  follows from a direct computation.

Next we prove  that for any given $\xi \in H^-$ the integral
equation (\ref{eq(3.2)}) has a unique solution in $C_{\eta}^+$. To
see this, let $J^s (u,\xi )$ denote  the right hand side of
equality (\ref{eq(3.2)}).  It is easy to see that $ J^s$ is
well-defined from $C_\eta^+ \times H^-$ to $C_\eta^+$. For each
$u, \bar u \in C_\eta^+$, we have  that
\begin{align}\label{eq(3.6)}
\begin{split}
&| J^s(u,\xi ) -  J^s (\bar u, \xi ) |_{C_\eta^+} \\
&\leq\sup_{t\in [0, \infty)} \Big\{e^{-\eta t - \int_0^t
z(\theta_s\omega) ds} \big(|\int_0^t
e^{A(t-s)+\int_s^tz(\theta_r\omega)dr}
P^-(G(\theta_s\omega,u)\\&\qquad\qquad\qquad\qquad\qquad\qquad\qquad\qquad\qquad- G(\theta_s\omega,\bar u))ds\\
&+\int^t_\infty e^{A(t-s)+\int_s^tz(\theta_r\omega)dr}
P^+(G(\theta_s\omega,u)-G(\theta_s\omega,\bar u))ds| \big)
\Big\}\\
&\leq\sup_{t\in [0, \infty)} \Big\{K \hbox{{\rm Lip}}_u G |u-\bar
u|_{C_\eta^+}\big(\int_0^t e^{(\beta-\eta)(t-s)}ds+\int^t_\infty
e^{(\alpha-\eta)(t-s)}ds\big)
\Big\}\\
&\leq K \hbox{{\rm Lip}}_u
G(\frac{1}{\eta-\beta}+\frac{1}{\alpha-\eta}) |u-\bar
u|_{C_\eta^+}.
\end{split}
\end{align}

Obviously $ J^s$ is Lipschitz continuous in $\xi$.  By the
assumption,
 $K \hbox{{\rm Lip}}_u G(\frac{1}{\eta-\beta}+\frac{1}{\alpha-\eta})< 1$,
 hence $ J^s$ is a uniform contraction with respect to the parameter
$\xi$.  By the uniform contraction mapping principle, we have that
for each $\xi \in H^-$,  the mapping $J^s (\cdot , \xi )$ has a
unique fixed point $u(\cdot; \xi, \omega ) \in C_\eta^+$ and
$u(\cdot ;\cdot, \omega)$ is Lipschitz from $H^-$ to $C_\eta^+$,
that is, $u(\cdot;\cdot,\omega )\in C_\eta^+$ is a unique solution
of the integral equation (\ref{eq(3.2)}). Furthermore one has for
the fixed point $u$ the estimate
\begin{align}\label{eq(3.7)}
|u(\cdot ;\xi, \omega) - u (\cdot ; \bar \xi, \omega
)|_{C_\eta^+}\leq \frac{K}{1-K \hbox{{\rm Lip}}_u
G(\frac{1}{\eta-\beta}+\frac{1}{\alpha-\eta})}|\xi - \bar\xi |.
\end{align}

Since $u(\cdot ;\xi, \omega)$ can be an $\omega$-wise limit of the
iteration of contraction mapping $J^s$ starting at $0$ and $J^s$
maps a $\mathcal{F}$-measurable function to a measurable function,
$u(\cdot ;\xi, \omega)$ is $\mathcal{F}$-measurable. On the other
hand, since $u(\cdot ;\xi, \omega)$ is Lipschitz continuous, by
Castaing and Valadier \cite{CasVal77}, Lemma III.14, the above
terms are measurable with respect to $(\xi, \omega,y)$.

Let $h^s (\xi,\omega) = P^+ u (0; \xi,\omega)$. Then
\[
h^s(\xi, \omega) =\int^0_\infty e^{-As\int_s^0z(\theta_r\omega)dr}
P^+G(\theta_s\omega,u(s;\xi,\omega))ds
\]
and $h^s(0, \omega)=0$.

Thus, by using (\ref{eq(3.7)}), we obtain that
\[
|h^s(\xi,\omega) - h^s(\bar \xi, \omega)|\leq \frac{K^2\hbox{{\rm
Lip}}_u G}{(\alpha-\eta)\big(1-K \hbox{{\rm Lip}}_u
G(\frac{1}{\eta-\beta}+\frac{1}{\alpha-\eta}\big)}|\xi - \bar\xi |
\]
and $h^s$ is measurable.  From the definition of $h^s(\xi,
\omega)$ and the claim that $u^0 \in M^s(\omega)$  if and only if
there exists $u \in C_\eta^+$ with $u(0)=u_0$ and satisfies
(\ref{eq(3.2)}) it follows that $u^0 \in M^s(\omega)$ if and only
if there exists $\xi \in H^-$ such that $u^0=\xi+h^s(\xi,
\omega)$, therefore,
$$
M^s(\omega) = \{\xi+ h^s(\xi, \omega) |\xi \in H^-\}.
$$

In order to see that $M^s(\omega)$ is a random set we need to show
that for any $x\in H$
\begin{equation}\label{eq(3.8)}
\omega\to\inf_{y\in H}\|x-(P^-y+ h^s(P^-y, \omega))\|
\end{equation}
is measurable, see Castaing and Valadier \cite{CasVal77}, Theorem
III.9. Let $H_c$ be a countable dense set of the separable space
$H$. Then the right hand side of (\ref{eq(3.8)}) is equal to
\begin{equation}\label{eq(3.9)}
\inf_{y\in H_c}\|x-P^-y+ h^s(P^-y, \omega))|
\end{equation}
which follows immediately by the continuity of $h^s(\cdot,
\omega)$. The measurability of any expression under the infimum of
(\ref{eq(3.8)}) follows since
$\omega\to h^s(P^-y,\omega)$ is measurable for any $y\in H$. \\

Finally, we show that $M^s(\omega)$ is invariant, i.e., for each
$u_0\in M^s(\omega)$, $u(s, u_0, \omega) \in M^s(\theta_s\omega)$
for all $ s\geq 0$. We first note that for each fixed $s \geq 0$,
$u(t+s, u_0,\omega)$ is a solution of
\[
\frac{du}{dt}=Au+z(\theta_t(\theta_s\omega))u+G(\theta_t(\theta_s\omega),u),\quad
u(0)=u(s, u_0, \omega).
\]
Thus, $u(t,u(s,u_0,\omega), \theta_s\omega)=u(t+s, u_0,\omega)$.

Since $u(\cdot, u_0, \omega) \in C_\eta^+$, $u(t,u(s,u_0,\omega),
\theta_s\omega) \in C_\eta^+$. Therefore, $u(s, u_0, \omega) \in
M^s(\theta_s\omega)$

This completes the proof.
\end{proof}

\begin{theorem} \label{Thm(3.2)}

$\tilde M^s(\omega)=T^{-1}(\omega,M^s(\omega))$ is a Lipschiz
stable manifold of the stochastic partial differential equation
(\ref{eq(2.2)}).
\end{theorem}
\begin{proof}
Let $u(t,\omega, x)$ denote the solution of  (\ref{eq(2.5)}) and
$\tilde u(t, \omega, x)$ denote the solution of (\ref{eq(2.2)}).
From Lemma \ref{lem(2.2)}, we have
\begin{align*}
\tilde u(t,\omega,&\tilde  M^s(\omega))=T^{-1}(\theta_t\omega,u(t,\omega,T(\omega,\tilde M^s(\omega))))\\
&=T^{-1}(\theta_t\omega, u(t,\omega,M^s(\omega)))\subset
T^{-1}(\theta_t\omega,M^s(\theta_t\omega)) =\tilde
M^s(\theta_t\omega).
\end{align*}
Hence, $\tilde M^s(\omega)$ is an invariant set.  We also notice
that
\begin{align*}
&\tilde M^s(\omega)\\ &= T^{-1}(\omega,M^s(\omega))\\ &
=\big\{u_0=T^{-1}(\omega, \xi+h^s(\xi, \omega)\big |\; \xi \in H^-
\big\} \\
&=\big\{u_0= e^{z(\omega)}(\xi+h^s(\xi, \omega))\big |\; \xi \in
H^-
\big\}\\
&=\big\{u_0= (\xi+h^s(e^{-z(\omega)}\xi, \omega))\big |\; \xi \in
H^- \big\}
\end{align*}
which implies that $\tilde M^s(\omega)$ is a Lipschitz stable
manifold given by the graph of a Lipschitz continuous function
$\tilde h^s(\xi, \omega)=h^s(e^{-z(\omega)}\xi, \omega)$ over the
space $H^-$.
\end{proof}

%%%% Section 4 %%%%%
%%%%%%%%%%%%%%%%%%%%%%%%%%%%%
\section{\bf Smoothness of Stable Manifolds}

In this section, we prove that for each $\omega$, $M^s(\omega)$ is
a $C^k$ smooth manifold. We have

\begin{theorem} \label{Thm(4.1)}   Assume that $G$ is $C^k$ in $u$. If $\beta < k\eta < \alpha$ and
\[
K \hbox{{\rm Lip}}_u
G(\frac{1}{i\eta-\beta}+\frac{1}{\alpha-i\eta})< 1 \quad\hbox{for
all}\; 1\leq i\leq k,
\]
then $M^s(\omega)$ is a $C^k$ invariant stable manifold for the
random
 partial differential  equation (\ref{eq(3.1)}), i.e., $h(\xi, \omega)$ is
$C^k$ in $\xi$.
\end{theorem}

\begin{proof} We prove this theorem by induction.
First, we consider $k=1$. Since
\[
K \hbox{{\rm Lip}}_u
G(\frac{1}{\eta-\beta}+\frac{1}{\alpha-\eta})< 1
\]
there exists a small number $\delta>0$ such that $\beta <
\eta-2\delta$ and
\[
K \hbox{{\rm Lip}}_u
G(\frac{1}{(\eta-\gamma)-\beta}+\frac{1}{\alpha-(\eta-\gamma)})<
1\quad \hbox{for all } 0 \leq \gamma \leq 2\delta.
\]
Thus, $J^s(\cdot,\xi, \omega)$ defined in the proof of Theorem
\ref{Thm(3.1)} is a uniform contraction in $C_{\eta-\gamma}^+
\subset  C_\eta^+$ for any $0\leq \gamma\leq 2\delta$. Therefore,
$u(\cdot;\xi,\omega)\in C_{\eta-\gamma}^+$. For $\xi_0\in H^-$, we
set
\[
S=e^{At+\int_0^tz(\theta_s\omega)ds},
\]
and
\begin{align*}
T v&=\int_0^t e^{A(t-s)+\int_s^tz(\theta_r\omega)dr}
P^-D_uG(\theta_s\omega,u(s;\xi_0,\omega)) v ds \\
&+\int_\infty^t e^{A(t-s)+\int_s^tz(\theta_r\omega)dr}
P^+D_uG(\theta_s\omega,u(s;\xi_0,\omega)) v ds
\end{align*}
for $v\in C_{\eta-\delta}^+$. From the assumption, we have that
$S$ is a bounded linear operator from $H^-$ to
$C_{\eta-\delta}^+$. Using the same arguments as we proved that
$J^s$ is a contraction, we have that $T$ is a bounded linear
operator from $C_{\eta-\delta}^+$ to itself and
\[
||T|| \leq K \hbox{{\rm Lip}}_u
G(\frac{1}{(\eta-\delta)-\beta}+\frac{1}{\alpha-(\eta-\delta)})<
1,
\]
which implies that $Id-T$ is invertible in $C_{\eta-\delta}^+$.
For $\xi, \xi_0\in H^-$, we set
\begin{align*}
I=&\int_0^t e^{A(t-s)+\int_s^tz(\theta_r\omega)dr}
P^-\Big[G(\theta_s\omega,u(s;\xi,\omega))-
G(\theta_s\omega,u(s;\xi_0,\omega))\\
&\quad\quad\quad\quad\quad\quad-D_uG(\theta_s\omega,u(s;\xi_0,\omega))
(u(s;\xi,\omega)-u(s;\xi_0,\omega))\Big]ds\\
&+\int^t_\infty e^{A(t-s)+\int_s^tz(\theta_r\omega)dr}
P^+\Big[G(\theta_s\omega,u(s;\xi,\omega))-
G(\theta_s\omega,u(s;\xi_0,\omega))\\
&\quad\quad\quad\quad\quad\quad-D_uG(\theta_s\omega,u(s;\xi_0,\omega))
(u(s;\xi,\omega)-u(s;\xi_0,\omega))\Big]ds.
\end{align*}

We claim that $|I|_{C_{\eta-\delta}^+}=o(|\xi-\xi_0|)$ as $\xi\to
\xi_0$. Using this claim, we obtain
\begin{align}\label{eq(4.1)}
\begin{split}
&u(\cdot;\xi, \omega)-u(\cdot;\xi_0, \omega)-T (u(\cdot;\xi,
\omega)-u(\cdot;\xi_0, \omega))
\\
&=S(\xi-\xi_0) + I \\
&=S(\xi-\xi_0) + o(|\xi-\xi_0|), \,\hbox{as}\,\xi \to \xi_0.\\
\end{split}
\end{align}
which yields
\begin{align*}
u(\cdot;\xi, \omega)-u(\cdot;\xi_0,
\omega)=(Id-T)^{-1}S(\xi-\xi_0)+o(|\xi-\xi_0|).
\end{align*}
Hence,  $u(\cdot;\xi, \omega) $ is differentiable in $\xi$ and its
derivative satisfies $D_\xi u(t;\xi,\omega)\in L(H^-,
C_{\eta-\delta}^+)$,  where $L(H^-, C_{\eta-\delta}^+)$ is the
usual space of bounded linear operators and

\begin{align}\label{eq(4.2)}
\begin{split}
D_\xi u(t;\xi,\omega) &=e^{At+\int_0^tz(\theta_s\omega)ds}P^-\cdot \\
& \quad +\int_0^t e^{A(t-s)\int_s^t+z(\theta_r\omega)dr}
P^-D_uG(\theta_s\omega,u(s;\xi,\omega))D_\xi u(s;\xi,\omega)ds\\
&\quad +\int_\infty^t e^{A(t-s)+\int_s^tz(\theta_r\omega)dr}
P^+D_uG(\theta_s\omega,u(s;\xi,\omega)) D_\xi u(s;\xi,\omega) ds
\end{split}
\end{align}

Now we prove that $|I|_{C_{\eta-\delta}^+} = o(|\xi-\xi_0|)$ as
$\xi\to \xi_0$. Let $N$ be a large  positive number to be chosen
later and let
\begin{align*}
I_1=&e^{-(\eta-\delta) t-\int_0^t z(\theta_s\omega) ds}
\Big\{|\int_N^t e^{A(t-s)+\int_s^tz(\theta_r\omega)dr}
P^-\Big[G(\theta_s\omega,u(s;\xi,\omega)) \\ &-
G(\theta_s\omega,u(s;\xi_0,\omega))-D_uG(\theta_s\omega,u(s;\xi_0,\omega))
(u(s;\xi,\omega)-u(s;\xi_0,\omega))\Big]ds|\Big\}\\
\end{align*}
for $t\geq N$ and $I_1=0$ for $t < N$;
\begin{align*}
I_2=&e^{-(\eta-\delta) t-\int^t_0 z(\theta_s\omega) ds}
\Big\{|\int_0^N e^{A(t-s)+\int_s^tz(\theta_r\omega)dr}
P^-\Big[G(\theta_s\omega,u(s;\xi,\omega)) \\ &-
G(\theta_s\omega,u(s;\xi_0,\omega))-D_uG(\theta_s\omega,u(s;\xi_0,\omega))
(u(s;\xi,\omega)-u(s;\xi_0,\omega))\Big]ds|\Big\}.\\
\end{align*}
Let $\bar N$ be a large  positive number to be chosen later. For
$0\leq t \leq \bar N$, we set
\begin{align*}
I_3=&e^{-(\eta-\delta) t-\int_t^0 z(\theta_s\omega) ds}
\Big\{|\int_{\bar N}^t e^{A(t-s)+\int_s^tz(\theta_r\omega)dr}
P^+\Big[G(\theta_s\omega,u(s;\xi,\omega))\\ &-
G(\theta_s\omega,u(s;\xi_0,\omega))-D_uG(\theta_s\omega,u(s;\xi_0,\omega))
(u(s;\xi,\omega)-u(s;\xi_0,\omega))\Big]ds|\Big\};\\
\end{align*}
\begin{align*}
I_4=&e^{-(\eta-\delta) t-\int_t^0 z(\theta_s\omega) ds}
\Big\{|\int_\infty^{\bar N} e^{A(t-s)+\int_s^tz(\theta_r\omega)dr}
P^+\Big[G(\theta_s\omega,u(s;\xi,\omega)) \\ &-
G(\theta_s\omega,u(s;\xi_0,\omega))-D_uG(\theta_s\omega,u(s;\xi_0,\omega))
(u(s;\xi,\omega)-u(s;\xi_0,\omega))\Big]ds|\Big\}.\\
\end{align*}
For $t \geq \bar N$, we set
\begin{align*}
I_5=&e^{-(\eta-\delta) t-\int_t^0 z(\theta_s\omega) ds}
\Big\{|\int_\infty^t e^{A(t-s)+\int_s^tz(\theta_r\omega)dr}
P^+\Big[G(\theta_s\omega,u(s;\xi,\omega)) \\ &-
G(\theta_s\omega,u(s;\xi_0,\omega))-D_uG(\theta_s\omega,u(s;\xi_0,\omega))
(u(s;\xi,\omega)-u(s;\xi_0,\omega))\Big]ds|\Big\}.\\
\end{align*}

It is sufficient to show that for any $\epsilon >0$ there is a
$\sigma
>0$ such that if $ |\xi-\xi_0| \leq \sigma $, then
$|I|_{C_{\eta-\delta}^+} \leq \epsilon |\xi-\xi_0|$. Note that

$$
|I|_{C_{\eta-\delta}^+} \leq \sup_{t\geq 0} I_1 + \sup_{t\geq 0}
I_2 + \sup_{0\leq t\leq \bar N} I_3+\sup_{0\leq t\leq \bar N} I_4
+ \sup_{t \geq \bar N} I_5.
$$
A computation similar to (\ref{eq(3.7)}) implies that

\begin{align*}
I_1 &\leq 2K\hbox{{\rm Lip}}_uG \int_N^t
e^{(\beta-(\eta-\delta))(t-s)}e^{-\delta s}
|u(\cdot;\xi,\omega)-u(\cdot;\xi_0,\omega)|_{C_{\eta-2\delta}^+} ds \\
& \leq \frac{2K^2\hbox{{\rm Lip}}_uG e^{-\delta
N}}{(\eta-\delta-\beta)(1-K \hbox{{\rm Lip}}_u
G(\frac{1}{\eta-2\delta-\beta}+\frac{1}{\alpha-(\eta-2\delta)}))}
|\xi - \xi_0 |.
\end{align*}

 Choose
$N$ so large that
\begin{align*}
\frac{2K^2\hbox{{\rm Lip}}_uG e^{-\delta
N}}{(\eta-\delta-\beta)(1-K \hbox{{\rm Lip}}_u
G(\frac{1}{\eta-2\delta-\beta}+\frac{1}{\alpha-(\eta-2\delta)}))}
\leq \frac{1}{4} \epsilon.
\end{align*}

Hence for such $N$ we have that
$$
\sup_{t\geq 0}I_1 \leq \frac{1}{4}\epsilon |\xi-\xi_0|_X.
$$
Fixing such $N$, for $I_2$ we have that
\begin{align*}
I_2 &\leq K\int_0^N e^{(\beta-(\eta-\delta))(t-s)}\Big\{
\int^1_0\big[|
D_uG(\theta_s\omega,\tau u(s;\xi,\omega)+(1-\tau)u(s;\xi_0,\omega)) \\
&-D_uG(\theta_s\omega,u(s;\xi_0,\omega))|\big] d\tau\Big\}
|u(\cdot;\xi,\omega)-u(\cdot;\xi_0,\omega)|_{C_{\eta-\delta}^+} ds \\
& \leq \frac{K^2|\xi-\xi_0|}{1-K \hbox{{\rm Lip}}_u
G(\frac{1}{\eta-\delta-\beta}+\frac{1}{\alpha-(\eta-\delta)})}\\&
\int_0^N e^{-(\beta-(\eta-\delta))s}\Big\{ \int^1_0\big[|
D_uG(\theta_s\omega,\tau
u(s;\xi,\omega)+(1-\tau)u(s;\xi_0,\omega))\\&
-D_uG(\theta_s\omega,u(s;\xi_0,\omega))|\big] d\tau\Big\}ds.\;  \\
\end{align*}

The last integral is on the compact interval $[0,N]$. Thus, from
the continuity of  the integrand  $(s, \xi)$,   we have that there
is a $\sigma_1
>0$ such that if $|\xi-\xi_0| \leq \sigma_1$, then
$$
\sup_{t\geq 0}I_2 \leq \frac{1}{4}\epsilon |\xi-\xi_0|.
$$
Therefore, if  $|\xi-\xi_0| \leq \sigma_1$, then
$$
\sup_{t\geq 0}I_1 + \sup_{t\geq 0}I_2 \leq \frac{1}{2}\epsilon
|\xi-\xi_0|.
$$
Similarly, by choosing $\bar N$ to be sufficiently large, we have
$$
\sup_{0\leq t\leq \bar N} I_4 + \sup_{t\geq \bar N}I_5 \leq
\frac{1}{4}\epsilon|\xi-\xi_0|,
$$
and for fixed such $\bar N$, there exists $\sigma_2>0$ such that
if $|\xi-\xi_0| \leq \sigma_2$, then
$$
\sup_{0\leq t\leq \bar N}I_3  \leq \frac{1}{4}\epsilon
|\xi_1-\xi_2|.
$$
Taking $\sigma = \min \{\sigma_1, \sigma_2\}$, we have that if
$|\xi-\xi_0| \leq \sigma$, then
$$
|I|_{C_{\eta-\delta}^+} \leq \epsilon |\xi-\xi_0|.
$$
Therefore  $|I|_{C_{\eta-\delta}^+} = o(|\xi-\xi_0|)$ as $\xi\to
\xi_0$. We now prove that $D_\xi u(\cdot;\cdot,\omega)$ is
continuous from $H^-$ to $C_\eta^+)$. For $\xi,\; \xi_0 \in H^-$,
using (\ref{eq(4.2)}), we have

\begin{align}\label{eq(4.3)}
\begin{split}
&D_\xi u(t;\xi,\omega)-D_\xi u(t;\xi_0,\omega)\\
&=\int_0^t e^{A(t-s)+\int_s^tz(\theta_r\omega)dr}
P^-\Big(D_uG(\theta_s\omega,u(s;\xi,\omega))D_\xi
u(s;\xi,\omega)\\&\hskip1.5in -
D_uG(\theta_s\omega,u(s;\xi_0,\omega))D_\xi u(s;\xi_0,\omega)\Big)ds\\
&\quad+\int_\infty^t e^{A(t-s)+\int_s^tz(\theta_r\omega)dr} P^+
\Big(D_uG(\theta_s\omega,u(s;\xi,\omega))D_\xi
u(s;\xi,\omega)\\&\hskip1.5in -
D_uG(\theta_s\omega,u(s;\xi_0,\omega))D_\xi
u(s;\xi_0,\omega)\Big)ds \\
&=\int_0^te^{A(t-s)+\int_s^tz(\theta_r\omega)dr}
P^-\Big(D_uG(\theta_s\omega,u(s;\xi,\omega))
\\&\hskip1.5in(D_\xi
u(s;\xi,\omega)- D_\xi u(s;\xi_0,\omega))\Big)ds\\
&\quad+\int_\infty^t e^{A(t-s)+\int_s^tz(\theta_r\omega)dr} P^+
\Big(D_uG(\theta_s\omega,u(s;\xi,\omega))
\\&\hskip1.5in (D_\xi u(s;\xi,\omega)-
D_\xi u(s;\xi_0,\omega))\Big)ds + \bar I,
\end{split}
\end{align}
where
\begin{align*} \bar I =&
\int_0^t e^{A(t-s)+\int_s^tz(\theta_r\omega)dr}
P^-\big(D_uG(\theta_s\omega,u(s;\xi,\omega))\\ & \hskip1.5in
-D_uG(\theta_s\omega,u(s;\xi_0,\omega)
\big)D_\xi u(s;\xi_0,\omega)ds\\
& +\int_\infty^t e^{A(t-s)+\int_s^tz(\theta_r\omega)dr}
P^+\big(D_uG(\theta_s\omega,u(s;\xi,\omega))\\ & \hskip1.5in
-D_uG(\theta_s\omega,u(s;\xi_0,\omega)
\big)D_\xi u(s;\xi_0,\omega)ds.\\
\end{align*}
 Then, estimating $|D_\xi u(\cdot;\xi,\omega)-D_\xi u(\cdot;\xi_0,\omega)|_{L(H^-,
 C_\eta^+)}$, we have
\begin{align*}
&|D_\xi u(\cdot;\xi,\omega)-D_\xi u(\cdot;\xi_0,\omega)|_{L(H^-, C_\eta^+)} \\
& \leq \frac{|\bar I|_{L(H^-, C_\eta^+)}}{1- K \hbox{{\rm Lip}}_u
G(\frac{1}{\eta-\beta}+\frac{1}{\alpha-\eta})}.
\end{align*}

Using the same argument we used for the last claim, we obtain that
$|\bar I|_{L(H^-, C_\eta^+)} = o(1)$ as $\xi \to \xi_0$. Hence
$D_\xi u(\cdot;\cdot,\omega)$ is continuous from $H^-$ to ${L(H^-,
C_\eta^+)}$. Therefore, $u(\cdot;\cdot,\omega)$ is $C^1$ from
$H^-$ to $C_\eta^+$. Now we show that $u$ is $C^k$ from $H^-$ to
$C^+_{k\eta}$ by induction for $k\ge 2$. By the induction
assumption, we know that $u$ is $C^{k-1}$ from $H^-$ to
$C^+_{(k-1)\eta}$ and the  and $(k-1)$-derivative $D^{k-1}_\xi
u(t;\xi,\omega)$ satisfies the following equation
\begin{align*}
D_\xi^{k-1} u =& \int^t_0e^{A(t-s)\int_s^t+z(\theta_r\omega)dr}
P^-(D_uG(\theta_s\omega,u)D_\xi^{k-1} u ds\\
&+ \int^t_\infty e^{A(t-s)+\int_s^tz(\theta_r\omega)dr}
P^+D_uG(\theta_s\omega,u)D_\xi^{k-1} u ds\\
& \int^t_0e^{A(t-s)\int_s^t+z(\theta_r\omega)dr}
P^-R_{k-1}(s, \xi, \omega) ds\\
&+ \int^t_\infty e^{A(t-s)+\int_s^tz(\theta_r\omega)dr}
P^+ R_{k-1}(s, \xi, \omega) ds\\
\end{align*}
where
$$
R_{k-1}(s,\xi,\omega) =  \sum_{i=0}^{k-3}\begin{pmatrix} k-2\\i
\end{pmatrix} D_\xi^{k-2-i} \big(D_uG(\theta_s
\omega,u(s;\xi,\omega))\big) D_\xi^{i+1}u(s;\xi,\omega).
$$
We note that $D^i_\xi u \in C^+_{i\eta} $ for $i=1,\cdots, k-1$
from the induction hypothesis. Thus, using the fact that $G$ is
$C^k$, we can verify that $R_{k-1}(\cdot,\xi,\omega) \in
L^{k-1}\big(H^-,C^+_{(k-1)\eta}\big)$ and is $C^1$ in $\xi$, where
$L^{k-1}\big(H^-,C^+_{(k-1)\eta}\big)$ is the usual space of
bounded $k-1$ linear forms. In order to insure that the above
integrals are well-defined one has to require that $
\beta<(k-1)\eta<\alpha$. This is the reason why we need the gap
condition. The fact that $t\to z(\theta_t\omega)$ has a sublinear
growth rate is also used in these analysis.
 Note
that from the assumption $\beta < k\eta < \alpha$ and
\[
K \hbox{{\rm Lip}}_u
G(\frac{1}{i\eta-\beta}+\frac{1}{\alpha-i\eta})< 1 \quad\hbox{for
all}\; 1\leq i\leq k.
\]
Using this fact and the same argument which we used in the case
$k=1$, we can show that $ D^{k-1}_\xi u(\cdot;\cdot,\omega)$ is
$C^1$ from $X$ to $L^k(H^-,C^+_{k\eta})$.   This completes the
proof.
\end{proof}

\begin{theorem} \label{Thm(4.2)}
Assume that $F(u)$ is $C^k$  smooth.
If $\beta < k\eta < \alpha$ and
\[
K \hbox{{\rm Lip}}_u
G(\frac{1}{i\eta-\beta}+\frac{1}{\alpha-i\eta})< 1 \quad\hbox{for
all}\; 1\leq i\leq k,
\]
then
$\tilde M^s(\omega)=T^{-1}(\omega,M^s(\omega))$  is a $C^k$  invariant
stable manifold for the stochastic partial differential  equation  (\ref{eq(2.2)}).
\end{theorem}
\begin{proof}
Since
\[
\tilde M^s(\omega)=\big\{\xi+\tilde h^s(\xi, \omega) \big |\;
\xi\in H^-\big\},
\]
$\tilde h^s(\xi, \omega)=e^zh^s(e^{-z(\omega)}\xi, \omega)$, and
$h^s(\xi, \omega)$ is $C^k$ in $\xi$, $\tilde h^s(\xi, \omega)$ is
$C^k$ in $\xi$.
\end{proof}

%%%%%%%%%%%%%%%%
%%%% Section 5 %%%%%%
%%%%%%%%%%%%%%%%
\section{\bf Smooth Unstable Manifolds}

All results obtained in Section 3 and Section 4 also hold for
unstable manifolds.

\begin{theorem} \label{Thm(5.1)}    If $$K \hbox{{\rm Lip}}_u G(\frac{1}{\eta-\beta}+\frac{1}{\alpha-\eta})<
1,$$ then there exists a Lipschitz unstable  manifold for
 the random partial differential  equation (\ref{eq(3.1)}), which is given by
\[
 M^u(\omega) = \{\xi+h^u(\xi, \omega)\big | \xi \in H^+\},
\]
where $h^u : H^+\to  H^-$ is a Lipschitz continuous mapping and
satisfies $h^u(0)=0$.
Moreover,
$\tilde M^u(\omega)=T^{-1}(\omega,M^u(\omega))$ is a Lipschiz
stable manifold of the stochastic  partial differential     equation
 (\ref{eq(2.2)}).
\end{theorem}

\begin{theorem} \label{Thm(5.2)}   Assume that the nonlinear term
$F$ and thus
$G$ is $C^k$ in $u$. If $\beta < k\eta < \alpha$ and
\[
K \hbox{{\rm Lip}}_u
G(\frac{1}{i\eta-\beta}+\frac{1}{\alpha-i\eta})< 1 \quad\hbox{for
all}\; 1\leq i\leq k,
\]
then $M^u(\omega)$ is a $C^k$ unstable  manifold for
the random partial differential  equation (\ref{eq(3.1)}),
i.e., $h^u(\xi, \omega)$
is $C^k$ in $\xi$.
Moreover,
$\tilde M^u(\omega)=T^{-1}(\omega,M^u(\omega))$ is a $C^k$
unstable manifold for the stochastic  partial differential     equation
 (\ref{eq(2.2)}).
\end{theorem}

Generally, a few modifications are needed to adapt the proofs
presented in Section 3 and Section 4 to the case of unstable
manifold. The most significant differences are the integral
equation  (\ref{eq(3.2)}) and the associated function space. We shall outline
the proofs and leave the details to the interested reader.\\

Corresponding to space $C^+_\eta$,   we define the Banach Space
for each $\beta<\eta<\alpha$
\[
C_\eta^-=\{\phi:(-\infty, 0]\to H \;|\;\phi \hbox{  is continuous
and }\sup_{t\leq 0}e^{-\eta t - \int_0^t z(\theta_\tau\omega)
d\tau} |\phi(t)| < \infty\}
\]
with the norm
\[|\phi|_{C_\eta^-}=\sup_{t\leq 0}e^{-\eta t-\int_0^t z(\theta_\tau\omega)
d\tau} |\phi(t)|.
\]

Let
\[M^u(\omega)=\{u_0\in H\;|\; u(\cdot, \omega, u_0) \in C_\eta^-\}
\]
This is the set of all initial datum through  which solutions
decay  as $e^{\eta t +\int_0^t z(\theta_\tau\omega) d\tau}$ as
$t\to -\infty$.\\

Clearly, $M^u(\omega)$ is invariant. In order to show that
$M^u(\omega)$ is given by the graph of a $C^k$ (or Lipschitz)
function, one needs to prove that $u^0 \in \ M^u(\omega)$ if and
only if there exists a function $ u(\cdot)\in C_\eta^-$ with $u(0)
= u^0$ and satisfies
\begin{align}\label{eq(5.1)}\begin{split}
u(t)=&e^{At+\int_0^tz(\theta_s\omega)ds}\xi+\int_0^t
e^{A(t-s)+\int_s^tz(\theta_r\omega)dr}
P^+G(\theta_s\omega,u(s))ds\\
&+\int^t_{-\infty} e^{A(t-s)+\int_s^tz(\theta_r\omega)dr}
P^-G(\theta_s\omega,u(s))ds.
\end{split}
\end{align}
where $\xi = P^+ u^0$.\\

The next step is to show  that for any given $\xi \in H^+$ the
integral equation (\ref{eq(5.1)}) has a unique solution in
$C_{\eta}^-$. To see this, letting $J^u (u,\xi )$ denote  the
right hand side of integral equation  (\ref{eq(5.1)}),  one may show that
$J^u$ is a unform contraction. Hence, by the uniform contraction
mapping principle, we have that for each $\xi \in H^+$,  the
mapping $J^u (\cdot , \xi )$ has a unique fixed point $u(\cdot;
\xi, \omega ) \in C_\eta^-$ and $u(\cdot ;\cdot, \omega)$ is
Lipschitz from $H^+$ to $C_\eta^-$. Thus, $u(\cdot;\cdot,\omega
)\in C_\eta^-$ is a solution of integral equation (\ref{eq(5.1)}).\\

Let $h^u (\xi,\omega) = P^- u (0; \xi,\omega)$. Then
\[
h^u(\xi, \omega) =\int^0_{-\infty}
e^{-As\int_s^0+z(\theta_r\omega)dr}
P^+G(\theta_s\omega,u(s;\xi,\omega))ds
\]
and $h^u(0, \omega)=0$ if $F(0)=0$ or $G(\omega,0)=0$.

Therefore,
$$
 M^u(\omega) = \{\xi+ h^u(\xi, \omega) |\xi \in H^+\}.
$$

In the same fashion as the case for the smoothness of stable
manifold, one may show that $h^u$ is $C^k$ when the assumptions in
Theorem \ref{Thm(5.2)} hold.

%%%%%%%%%%%%%%%%%
%%%%%%%%%%%%%%%%%

\end{document}